\documentclass[12pt]{article}
\usepackage{graphicx}
\usepackage{amsmath,amsthm,amssymb,enumerate}%, esint}
\usepackage{euscript,mathrsfs}
\usepackage[left=2cm,right=2cm,top=3.5cm,bottom=3.5cm]{geometry}
\usepackage{color}
\catcode`\@=11 \@addtoreset{equation}{section}

\catcode`\@=12

\allowdisplaybreaks

\newtheorem{Theorem}{Theorem}[section]
\newtheorem{Proposition}[Theorem]{Proposition}
\newtheorem{Lemma}[Theorem]{Lemma}
\newtheorem{Corollary}[Theorem]{Corollary}

\theoremstyle{definition}
\newtheorem{Definition}[Theorem]{Definition}

\newtheorem{Remark}[Theorem]{Remark}

\newcommand{\bTheorem}[1]{
\begin{Theorem} \label{T#1} }
\newcommand{\eT}{\end{Theorem}}

\newcommand{\bProposition}[1]{
\begin{Proposition} \label{P#1}}
\newcommand{\eP}{\end{Proposition}}

\newcommand{\bLemma}[1]{
\begin{Lemma} \label{L#1} }
\newcommand{\eL}{\end{Lemma}}

\newcommand{\bCorollary}[1]{
\begin{Corollary} \label{C#1} }
\newcommand{\eC}{\end{Corollary}}

\newcommand{\bRemark}[1]{
\begin{Remark} \label{R#1} }
\newcommand{\eR}{\end{Remark}}

\newcommand{\bDefinition}[1]{
\begin{Definition} \label{D#1} }
\newcommand{\eD}{\end{Definition}}

\newcommand{\tvE}{\tilde{E}}

\newcommand{\tvm}{\tilde{\vc{m}}}
\newcommand{\bfphi}{\boldsymbol{\varphi}}

\newcommand{\bFormula}[1]{
\begin{equation} \label{#1}}
\newcommand{\eF}{\end{equation}}

\newcommand{\Ov}[1]{\overline{#1}}

\newcommand{\DC}{C^\infty_c}

\newcommand{\aleq}{\stackrel{<}{\sim}}

\newcommand{\vr}{\varrho}

\newcommand{\tvr}{\tilde \vr}
\newcommand{\tvu}{{\tilde \vu}}
\newcommand{\tvt}{\tilde \vt}

\newcommand{\vt}{\vartheta}
\newcommand{\vu}{\vc{u}}
\newcommand{\vm}{\vc{m}}

\newcommand{\vc}[1]{{\bf #1}}

\newcommand{\Div}{{\rm div}_x}
\newcommand{\Grad}{\nabla_x}

\newcommand{\dx}{\,{\rm d} {x}}

\newcommand{\dt}{\,{\rm d} t }

\newcommand{\intO}[1]{\int_{\Omega} #1 \ \dx}

\definecolor{Cgrey}{rgb}{0.85,0.85,0.85}
\definecolor{Cblue}{rgb}{0.50,0.85,0.85}
\definecolor{Cred}{rgb}{1,0,0}
\definecolor{fancy}{rgb}{0.10,0.85,0.10}

\newcommand\Cbox[2]{%
    \newbox\contentbox%
    \newbox\bkgdbox%
    \setbox\contentbox\hbox to \hsize{%
        \vtop{
            \kern\columnsep
            \hbox to \hsize{%
                \kern\columnsep%
                \advance\hsize by -2\columnsep%
                \setlength{\textwidth}{\hsize}%
                \vbox{
                    \parskip=\baselineskip
                    \parindent=0bp
                    #2
                }%
                \kern\columnsep%
            }%
            \kern\columnsep%
        }%
    }%
    \setbox\bkgdbox\vbox{
        \color{#1}
        \hrule width  \wd\contentbox %
               height \ht\contentbox %
               depth  \dp\contentbox
        \color{black}
    }%
    \wd\bkgdbox=0bp%
    \vbox{\hbox to \hsize{\box\bkgdbox\box\contentbox}}%
    \vskip\baselineskip%
}

\date{}

\begin{document}

%%%%%%%%%%%%%%%%%%%%%%%%%%%%%%%%

\title{Maximal dissipation principle for the complete Euler system}

\author{Jan B\v rezina \and Eduard Feireisl
\thanks{The research of E.F.~leading to these results has received funding from the
European Research Council under the European Union's Seventh
Framework Programme (FP7/2007-2013)/ ERC Grant Agreement
320078. The Institute of Mathematics of the Academy of Sciences of
the Czech Republic is supported by RVO:67985840.}
}

\date{\today}

\maketitle

\bigskip

\centerline{Tokyo Institute of Technology}

\centerline{ 2-12-1 Ookayama, Meguro-ku, Tokyo, 152-8550, Japan}

\bigskip

\centerline{Institute of Mathematics of the Academy of Sciences of the Czech Republic}

\centerline{\v Zitn\' a 25, CZ-115 67 Praha 1, Czech Republic}

\bigskip

\begin{abstract}

We introduce the concept of maximal dissipative measure--valued solution to the complete Euler system. These are solutions that maximize the entropy production rate. We show that these solutions exist under fairly general hypotheses imposed on the data and constitutive relations.

\end{abstract}


%\tableofcontents

\section{Introduction}
\label{i}

The principle of \emph{maximal dissipation} was proposed by Dafermos \cite{Dafer}, \cite{Daf55}, \cite{Daf4} as a natural admissibility criterion to rule out the nonphysical solutions to equations and systems in continuum dynamics. We examine this criterion in the context of the complete Euler system describing the motion of a compressible inviscid fluid:
\begin{equation}
\label{i1}
\begin{split}
\partial_t \vr + \Div (\vr \vu) &= 0,\\
\partial_t (\vr \vu) + \Div (\vr \vu \otimes \vu) + \Grad p(\vr, \vt) &= 0,\\
\partial_t \left( \frac{1}{2} \vr |\vu|^2 + \vr e(\vr, \vt) \right) +
\Div \left[ \left( \frac{1}{2} \vr |\vu|^2 + \vr e(\vr, \vt) + p(\vr, \vt) \right) \vu \right] &= 0.
\end{split}
\end{equation}
The unknowns are the mass density $\vr = \vr(t,x)$, the velocity field $\vu = \vu(t,x)$, and the (absolute) temperature $\vt = \vt(t,x)$, considered in the Eulerian reference frame $t = (0,T)$, $x \in \Omega \subset R^N$, $N=1,2,3$. The pressure $p$ and the specific internal energy $e$ are given functions of $(\vr,\vt)$ satisfying
Gibbs' equation
\begin{equation} \label{r3}
\vt D s = De + p D \left( \frac{1}{\vr} \right),
\end{equation}
where $s$ denotes the \emph{specific entropy}. In accordance with the Second law of thermodynamics, the entropy $s=s(\vr, \vt)$ satisfies the transport equation
\begin{equation} \label{i5}
\partial_t (\vr s(\vr, \vt)) + \Div (\vr s(\vr, \vt) \vu) = \sigma, \ \sigma \geq 0,
\end{equation}
where $\sigma$ is the entropy production rate. It is easy to check that $\sigma \equiv 0$ as long as all quantities in (\ref{i1}) are continuously differentiable.

As is well known, smooth solutions of (\ref{i1}) exist only for a finite lap of time after which singularities develop for a fairly generic class of initial data.
Therefore global-in-time solutions may exist only in a weak sense, where the derivatives in (\ref{i1}) are understood in the sense of distributions.
 An iconic
example is the Riemann problem completely understood in the 1-D setting, see e.g. Chen et al. \cite{CheFr2},
\cite{CheFr1}, where the physically relevant solutions satisfy the entropy balance (\ref{i5}) with $\sigma \ne 0$. Moreover, the unique solution can be singled out by maximizing the entropy production rate
$\sigma$, see Dafermos \cite{Dafer}.

The situation turns more complex in the multidimensional setting. As shown in \cite{FeKlKrMa}, problem (\ref{i1}) is basically ill--posed -
admits infinitely many weak solutions - in the class of bounded measurable functions for a large class of initial data. In addition, all of these solutions
satisfy (\ref{i5}) with $\sigma = 0$. This fact strongly suggests that maximizing $\sigma$ may rule out at least some of the possibly nonphysical solutions.

The entropy production rate $\sigma$ - a non--negative distribution - can be interpreted as a non--negative Borel measure sitting on the physical
space $[0,T] \times \Ov{\Omega}$. Let $(\vr_1, \vt_1, \vu_1)$, $(\vr_2, \vt_2, \vu_2)$ be two solutions of the Euler system (\ref{i1})
with the entropy production rates $\sigma_1$, $\sigma_2$, respectively. We say that
\begin{equation} \label{i6}
(\vr_1, \vt_1, \vu_1) \succsim (\vr_2, \vt_2, \vu_2) \ \mbox{iff} \ \ \sigma_1 \geq \sigma_2 \ \mbox{in}\ [0,T] \times \Ov{\Omega}.
\end{equation}
Our goal is to identify the class of \emph{maximal} solutions with respect to the relation $\succsim$. In contrast with the approach of Dafermos \cite{Dafer},
where maximality is understood \emph{globally} in the space variable, meaning maximizing $\intO{ \sigma }$ rather than $\sigma$, relation (\ref{i6})
requires maximality of $\sigma$ \emph{locally} on any subset of the associated physical space. Note that the ``global'' approach is probably too rough
to rule out the ``wild'' solutions, see Chiodaroli and Kreml \cite{ChiKre}.

Motivated by the pioneering paper of DiPerna \cite{DiP2}, we examine maximality of the entropy production rate in the class of \emph{dissipative measure--valued}
(DMV) \emph{solutions} introduced in \cite{BreFei17}. In particular, we show the \emph{existence} of a maximal DMV solution for a general class of initial data.
Similarly to \cite{DiP2}, the existence of a maximal solution is obtained by applying
Zorn's lemma (a variant of Axiom of Choice) argument. This is conditioned by uniform estimates for an ordered family of solutions and compactness in the class of DMV solutions. The paper is organized as follows. Preliminary material including a proper definition
of DMV solutions is collected in Section \ref{p}. Maximality is introduced and the main result stated in Section \ref{M}. Section \ref{E} contains the proof of the main result. Supplementary material is summarized in Section \ref{r} (Appendix).

\section{Preliminaries, dissipative measure--valued solutions}
\label{p}

For the sake of simplicity, we suppose that the pressure $p$ and the internal energy $e$ are interrelated through the caloric equation of state
\begin{equation} \label{r9}
p = (\gamma - 1) \vr e , \ \gamma > 1.
\end{equation}
More general equations of state can be treated in a similar fashion.

\subsection{Thermodynamic stability}

We impose the hypothesis of \emph{thermodynamic stability}:
\begin{equation} \label{r7}
\frac{\partial p}{\partial \vr}(\vr, \vt) > 0, \ \frac{\partial e}{\partial \vt} (\vr, \vt) > 0.
\end{equation}

In the context of weak or measure--valued solutions, it is convenient to consider the entropy $s = s(\vr, e)$ as
a function of $\vr$ and $e$. Accordingly, Gibbs' equation (\ref{r3}) yields
\begin{equation} \label{r5}
\frac{\partial s}{\partial e}(\vr, e) = \frac{1}{\vt}, \ \frac{\partial s}{\partial \vr}(\vr, e) = - \frac{p}{\vt \vr^2} = -
(\gamma - 1) \frac{e}{\vt \vr} = - (\gamma - 1) \frac{e}{\vr} \frac{\partial s}{\partial e}(\vr, e),
\end{equation}
where the first relation may be seen as the definition of the absolute temperature $\vt$. Moreover, the second relation in (\ref{r5})
represents a first order equation for $s$ that can be solved explicitly yielding
\begin{equation} \label{r10}
s(\vr, e) = S \left( \frac{(\gamma - 1) e}{\vr^{\gamma - 1}} \right) =
S \left( \frac{p}{\vr^{\gamma}} \right)
\end{equation}
for a certain function $S$. In accordance with (\ref{r5}), specifically with the requirement of positivity of the absolute temperature, we have
\[
S' > 0.
\]
Moreover, as shown in Appendix, the hypothesis of thermodynamic stability (\ref{r7}) implies
\[
\frac{p(\vr, \vt)}{\vr^\gamma} \searrow \Ov{p} \geq 0 \mbox{ as } \vt \to 0
\]
and
\[
(\gamma - 1) S'(Z) + \gamma S''(Z) Z < 0
\ \mbox{for all}\ Z > \Ov{p}.
\]
Finally, lifting $S$ by a constant as the case may be we may assume that
\[
\lim_{Z \to \Ov{p}+} S(Z) \in \left\{ 0, - \infty \right\}.
\]
Summarizing, we require the entropy $s$ to be given by (\ref{r10}), with
\begin{equation} \label{r11}
\begin{split}
S \in C^2(\Ov{p}, \infty),\ S'(Z) > 0&,\  (\gamma - 1) S'(Z) + \gamma S''(Z) Z < 0 \ \mbox{for}\ Z > \Ov{p},\\
\lim_{Z \to \Ov{p}+} S'(Z) = \infty&, \ \lim_{Z \to \Ov{p}+} S(Z) \in \{ 0, - \infty \}.
\end{split}
\end{equation}

\begin{Remark} \label{R1}

The standard example of $S$ corresponds to the perfect gas
\[
\Ov{p} = 0, \ S(Z) = \frac{1}{\gamma - 1} \log(Z),
\]
where $\lim_{Z \to 0+} S(Z) = - \infty$. The case $\lim_{Z \to \Ov{p}+} S(Z) = 0$ reflects the Third law of thermodynamics - the entropy vanishes
when the temperature approaches the absolute zero, cf. Belgiorno \cite{BEL1}, \cite{BEL2}. Finally, the case $\Ov{p} > 0$ corresponds to the
presence of ``cold'' pressure characteristic for the electron gas, see Ruggeri and Trovato \cite{RUTR}.

\end{Remark}

\subsection{Conservative variables}

To introduce the concept of measure--valued solution for complete Euler system, it is more convenient to formulate the problem in the \emph{conservative variables}:
\[
\vr, \ \vc{m}= \vr \vu, \ E = \frac{1}{2} \frac{|\vc{m}|^2}{\vr} + \vr e.
\]
The reason for changing the phase space is the fact that the temperature $\vt$ as well as the velocity 
$\vu$ may not be correctly defined on the (hypothetical) vacuum set.
As the measure--valued solutions are typically generated as weak limits of suitable approximation schemes, the presence of vacuum zones cannot be
{\it a priori} excluded.

The system (\ref{i1}) rewrites as
\begin{equation}
\label{p1}
\begin{split}
\partial_t \vr + \Div \vc{m} &= 0,\\
\partial_t \vc{m} + \Div \left( \frac{\vc{m} \otimes \vc{m}}{\vr} \right) + (\gamma - 1) \Grad \left( E - \frac{1}{2} \frac{|\vc{m}|^2}{\vr} \right)  &= 0,\\
\partial_t E +
\Div \left[ \left( E + (\gamma - 1) \left( E - \frac{1}{2} \frac{|\vc{m}|^2}{\vr} \right) \right) \frac{\vc{m}}{\vr} \right] &= 0,
\end{split}
\end{equation}
together with the associated entropy inequality
\begin{equation} \label{p2}
\partial_t \left( \vr S \left( (\gamma - 1) \frac{E - \frac{1}{2} \frac{|\vc{m}|^2}{\vr} }{\vr^\gamma} \right) \right)
+ \Div \left[ S \left( (\gamma - 1) \frac{E - \frac{1}{2} \frac{|\vc{m}|^2}{\vr} }{\vr^\gamma} \right) \vc{m} \right] \equiv \sigma \geq 0.
\end{equation}

Although the thermodynamic functions are well defined for regular values $\vr > 0$, $\vt > 0$ of the standard variables, where the latter condition corresponds in the conservative setting
to $E - \frac{1}{2} \frac{|\vc{m}|^2}{\vr} > 0$, we need them to be defined even for the limit values $\vr = 0$, $\vt = 0$. To that end, we first define the
kinetic energy,
\[
E_{\rm kin}(\vr, \vc{m}) = \frac{1}{2} \frac{|\vc{m}|^2}{\vr} = \left\{ \begin{array}{l} \frac{1}{2} \frac{|\vc{m}|^2 }{\vr} \ \mbox{for}\ \vr > 0,\\ \\
0 \ \mbox{if} \ \vc{m} = 0, \\ \\ \infty \ \mbox{otherwise.}
\end{array} \right.
\]
Note that $E_{\rm kin}$ is lower semi--continuous convex function defined on the set $\{ \vr \geq 0, \ \vc{m} \in R^N \}$. Similarly, we introduce the
total entropy
\[
\mathcal{S}(\vr, \vc{m}, E) = \vr S \left( (\gamma - 1) \frac{E - \frac{1}{2} \frac{|\vc{m}|^2}{\vr} }{\vr^\gamma} \right) = \vr S \left( (\gamma - 1) \frac{E - E_{\rm kin}(\vr, \vc{m}) }{\vr^\gamma} \right)
\]
\begin{equation} \label{p3}
\mathcal{S}(\vr, \vc{m}, E)
=
\left\{ \begin{array}{l}
\vr S \left( (\gamma - 1) \frac{E - \frac{1}{2} \frac{|\vc{m}|^2}{\vr} }{\vr^\gamma} \right) \ \mbox{if}\ \vr > 0,
\ E \geq \frac{1}{2} \frac{|\vc{m}|^2}{\vr} + \frac{ \Ov{p} }{\gamma - 1} \vr^\gamma,\\ \\
\lim_{\vr \to 0+} \vr S \left( (\gamma - 1) \frac{E}{\vr^\gamma } \right), \ \mbox{if}\ \vr = 0, \ \vc{m} = 0, \ E > 0,\\ \\
\lim_{E \to 0+} \left[ \lim_{\vr \to 0+} \vr S \left( (\gamma - 1) \frac{E}{\vr^\gamma } \right) \right]
\mbox{if} \ \vr = E = \vc{m} = 0, \\ \\
- \infty \ \mbox{otherwise}.
\end{array} \right.
\end{equation}
The total entropy $\mathcal{S}$ defined this way is a concave upper semi--continuous function defined on the set
$\{ \vr \geq 0,\ \vc{m} \in R^N, \ E \geq 0 \}$.

\subsection{DMV solutions}

The final objective of this preliminary part is to recall the definition of the dissipative measure--valued (DMV) solutions introduced in \cite{BreFei17} and \cite{BreFei17A}
respectively. The DMV solutions represent, roughly speaking, the most general object that complies with the principle of \emph{weak--strong} uniqueness. They coincide with a strong solution emanating from the same initial data as long as the latter exists, see \cite{BreFei17}. For the sake of simplicity, we consider the space periodic boundary conditions, meaning the physical domain $\Omega$ can be identified with the flat torus,
\[
\Omega = \left( [-1,1]|_{\{ - 1, 1\}} \right)^N.
\]
The definition can be easily adapted to the more realistic no--flux condition
\[
\vu \cdot \vc{n}|_{\partial \Omega} = 0
\]
or other types of admissible boundary behavior as the case may be.

The initial state of the system is given through a parameterized family of probability measures $\{ U_{0,x} \}_{x \in \Omega}$ defined on the
phase space
\[
\mathcal{Q} \equiv
\left\{ (\vr, \vc{m}, E) \ \Big| \ \vr \geq 0, \ \vc{m} \in R^N, \ E \geq 0 \right\},
\]
\[
(\vr_0, \vc{m}_0, E_0) (x) = \left[ \left< U_{0,x}; \vr \right>,
\left< U_{0,x}, \vc{m} \right>, \left< U_{0,x}; E \right> \right],
\]
where $\left< Y_{0,x}, g(\vr, \vc{m}, E) \right>$ denotes the expected value of a (Borel) function $g$ defined on $\mathcal{F}$. In addition, it is assumed
that the mapping $x \mapsto U_{0,x}$ belongs to $L^\infty_{{\rm weak-(*)}}(\Omega; \mathcal{P}(\mathcal{F}))$.

Similarly, a DMV solution is represented by a family of probability measures
\[
\{ U_{t,x} \}_{(t,x) \in (0,T) \times \Omega},\ U \in L^\infty_{{\rm weak}-(*)} ((0,T) \times \Omega; \mathcal{P}(\mathcal{F}) ),
\]
the non-linearities in (\ref{p1}), (\ref{p2}) are replaced by their expected values whereas the derivatives are understood in the sense of distributions.

\begin{Definition} \label{D1}
A parameterized family of probability measures $U \in L^\infty_{{\rm weak}-(*)} ((0,T) \times \Omega; \mathcal{P}(\mathcal{F}) )$ is
called a \emph{dissipative measure--valued (DMV) solution} to the Euler system (\ref{p1}), (\ref{p2}), with the initial data
$U_0 \in L^\infty_{{\rm weak}-(*)} (\Omega; \mathcal{P}(\mathcal{F}) )$ if the following holds:
\begin{itemize}
\item
\begin{equation} \label{p4}
\int_0^T \intO{ \left[ \left< U_{t,x}; \vr \right> \partial_t \varphi + \left< U_{t,x}; \vc{m} \right> \cdot
\Grad \varphi \right] } \dt = - \intO{ \left< U_{0,x} ; \vr \right> \varphi(0, \cdot)}
\end{equation}
for any $\varphi \in \DC([0,T) \times \Omega)$;
\item
\begin{equation} \label{p5}
\begin{split}
\int_0^T & \intO{ \left[ \left< U_{t,x}; \vc{m} \right> \cdot \partial_t \bfphi  + \left< U_{t,x}; \frac{ \vc{m} \otimes \vc{m} }{\vr}
\right> : \Grad \bfphi + (\gamma - 1) \left< U_{t,x}; E - \frac{1}{2} \frac{|\vc{m}|^2 }{\vr} \right> \Div \bfphi \right] }\dt \\
&= - \intO{ \left< U_{0,x}; \vc{m} \right> \cdot \bfphi(0, \cdot) } +
\int_0^T \int_\Omega \Grad \bfphi : {\rm d} \mu_C
\end{split}
\end{equation}
for any $\bfphi \in \DC([0, T) \times \Omega; R^N)$, where $\mu_C$ is a (vectorial) signed measure on $[0,T]\times \Omega$;
\item
\begin{equation} \label{p6}
\intO{ \left< U_{\tau,x}; E \right> } \leq \intO{ \left< U_{0,x}; E \right> } \ \mbox{for a.a.}\ \tau \in (0,T);
\end{equation}
\item
\begin{equation} \label{p7}
\begin{split}
\int_0^T &\intO{ \left[ \left< U_{t,x} ; \mathcal{S}(\vr, \vc{m}, E) \right> \partial_t \varphi +
\left< U_{t,x}; \mathcal{S} (\vr, \vc{m}, E) \frac{\vc{m}}{\vr} \right> \cdot \Grad \varphi  \right] } \dt \\
&\leq - \intO{ \left< U_{0,x}; \mathcal{S}(\vr, \vc{m}, E) \right> \varphi(0, \cdot) }
\end{split}
\end{equation}
for any $\varphi \in \DC([0,T) \times \Omega)$, $\varphi \geq 0$;
\item
\begin{equation} \label{p8}
\int_0^\tau \int_\Omega d \left| \mu_C \right|  \leq c(N, \gamma) \int_0^\tau \intO{ \left[ \left< U_{0,x}; E \right> - \left< U_{t,x}; E \right> \right]}
\dt \ \mbox{for any}\ 0 \leq \tau < T.
\end{equation}
\end{itemize}
\end{Definition}

Integrability of the convective term in the entropy inequality (\ref{p7}) may be sometimes problematic. In such a case, it is convenient to replace
the entropy $\mathcal{S}$ by its renormalization $\mathcal{S}_\chi$ defined as
\[
\mathcal{S}_\chi = \vr \chi \left( S \left( (\gamma - 1) \frac{E - \frac{1}{2} \frac{|\vc{m}|^2}{\vr} }{\vr^\gamma}          \right) \right),
\]
where $\chi$ is an increasing, concave function and $\chi (S) \leq \Ov{\chi}$ for all $S$. Note that composition
$\chi \circ S$ enjoys the same concavity properties as $S$ specified in (\ref{r11}), in particular, we may extend $\mathcal{S}_\chi$ to the full range
$\{ \vr \geq 0, \vc{m} \in R^N, E \geq 0 \}$ exactly as in (\ref{p3}). In addition, as $\chi$ is bounded from above, we get
\begin{equation} \label{p10}
\mathcal{S}_\chi = 0 \ \mbox{whenever}\ \vr = 0, \ \vc{m} = 0.
\end{equation}

Using $\mathcal{S}_\chi$ in place of $\mathcal{S}$ we may define renormalized dissipative measure--valued (rDMV) solutions as follows:
\begin{Definition} \label{D2}
A parameterized family of probability measures $U \in L^\infty_{{\rm weak}-(*)} ((0,T) \times \Omega; \mathcal{P}(\mathcal{F}) )$ is
called a \emph{renormalized dissipative measure--valued (rDMV) solution} to the Euler system (\ref{p1}), (\ref{p2}), with the initial data
$U_0 \in L^\infty_{{\rm weak}-(*)} (\Omega; \mathcal{P}(\mathcal{F}) )$ if it satisfies
the conditions (\ref{p4}--\ref{p6}), and (\ref{p8})
from Definition \ref{D1} whereas (\ref{p7}) is replaced by
\begin{equation} \label{p11}
\begin{split}
\int_0^T &\intO{ \left[ \left< U_{t,x} ; \mathcal{S}_\chi (\vr, \vc{m}, E) \right> \partial_t \varphi +
\left< U_{t,x}; \mathcal{S}_\chi (\vr, \vc{m}, E) \frac{\vc{m}}{\vr} \right> \cdot \Grad \varphi  \right] } \dt \\
&\leq - \intO{ \left< U_{0,x}; \mathcal{S}_\chi(\vr, \vc{m}, E) \right> \varphi(0, \cdot) }
\end{split}
\end{equation}
for any $\varphi \in \DC([0,T) \times \Omega)$, $\varphi \geq 0$, and any
\[
\chi \mbox{ defined on } ( S(\Ov{p}), \infty),  \ \mbox{increasing, concave},\ \chi(Z) \leq \Ov{\chi} \mbox{ for all } Z, \ \lim_{Z \to S(\Ov{p})+} \chi(Z) = S (\Ov{p}).
\]
\end{Definition}

\begin{Remark} \label{R}

In view of (\ref{p10}), we may write
\[
\mathcal{S}_\chi = \vr \chi \left( S \left( (\gamma - 1) \frac{E - \frac{1}{2} \frac{|\vc{m}|^2}{\vr} }{\vr^\gamma}  \right) \right)
= \vr \chi \left( S \left( (\gamma - 1) \frac{E - E_{\rm kin}(\vr, \vc{m}) }{\vr^\gamma}  \right) \right)
\]
in the full range $\vr \geq 0$, $\vc{m} \in R^N$, $E \geq \frac{1}{2} \frac{|\vc{m}|^2}{\vr} + \frac{\Ov{p}}{\gamma - 1} \vr^\gamma$.

\end{Remark}

Any ``standard'' weak solution $(\vr, \vc{m}, E)$ may be identified with a measure--valued solutions $U$ via
\[
U_{t,x} = \delta_{\vr(t,x), \vc{m}(t,x), E(t,x) } \ \mbox{for a.a.}\ (t,x) \in (0,T) \times \Omega,
\]
where $\delta_Z$ denotes the Dirac measure supported by $Z$.
As observed in \cite{BreFei17A}, the DMV solutions arise as zero dissipation limits of the Navier--Stokes--Fourier system, where renormalization of the
entropy equation is excluded for the primitive system, while rDMV solutions may be associated with the ``artificial'' viscosity limits for which
(\ref{p11}) holds. The natural inequality $E \geq \frac{1}{2} \frac{|\vc{m}|^2}{\vr} + \frac{\Ov{p}}{\gamma - 1} \vr^\gamma$ must be enforced through
the initial data, typically in a stronger form
\[
U_{0,x} \left\{  \vr > 0,\ (\gamma - 1) \frac{E - \frac{|\vc{m}|^2}{\vr} }{\vr^\gamma} \geq \Ov{p} + \delta \right\} = 1
\ \mbox{for a.a.}\ x \in \Omega, \ \delta > 0,
\]
which corresponds to positivity of the initial temperature. As shown in \cite{BreFei17}, the above property is propagated in time for rDMV solutions at least on the
non-vacuum set, meaning the conditional probability
\[
U_{t,x} \left\{ (\gamma - 1) \frac{E - \frac{|\vc{m}|^2}{\vr} }{\vr^\gamma} \geq \Ov{p} + \delta \ \Big| \vr > 0 \right\} = 1
\]
for a.a. $(t,x)$. In particular, we get
\[
U_{t,x} \left\{ \chi \circ S(\Ov{p} + \delta) \leq \chi\circ S \left( (\gamma - 1) \frac{E - \frac{1}{2} \frac{|\vc{m}|^2}{\vr} }{\vr^\gamma}  \right)
\leq \Ov{\chi} \right\} = 1 \ \mbox{for a.a.} \ (t,x) \in (0,T) \times \Omega.
\]

To conclude, we remark that the family of DMV and rDMV solutions for given initial data is closed with respect to convex combinations. In particular, in view of the results obtained in \cite{FeKlKrMa}, there is a vast class of initial data for which the Euler system admits infinitely many nontrivial DMV solutions. Here nontrivial means that  they do not consist of a single Dirac mass.

\subsection{Weak--strong uniqueness}

We conclude this preliminary part by reproducing the proof of the DMV - strong uniqueness principle from \cite{BreFei17}, here adapted to the
conservative variables $(\vr, \vc{m}, E)$.

\subsubsection{Relative energy}

We start by introducing the relative energy
\[
\mathcal{E} \left( \vr, \vc{m}, E \Big| \tilde \vr, \tilde \vt, \tilde \vu \right)  \equiv
E - \tilde \vt \mathcal{S}(\vr, \vc{m}, E) - \vc{m} \cdot \tilde \vu + \frac{1}{2} \vr |\tilde \vu|^2 +
p(\tilde \vr, \tilde \vt) - \left( e(\tilde \vr, \tilde \vt) - \tilde \vt s( \tilde \vr, \tilde \vt) + \frac{p(\tilde \vr, \tilde \vt)}{\tilde \vr} \right) \vr,
\]
cf. formula (\ref{r6}) in Appendix. Note that $\mathcal{E}$ should be seen as a function of six variables, namely $(\vr, \vc{m}, E, \tvr, \tvm, \tvE)$.
The DMV solutions satisfy the relative energy inequality:
\begin{equation} \label{rei}
\begin{split}
\Big[ \int_\Omega &\left< U_{t,x}; \mathcal{E} \left( \vr, \vc{m}, E \Big| \tilde \vr, \tilde \vt, \tilde \vu \right) \right> \ \dx \Big]_{t = 0}^{t = \tau}
+ \intO{ \left[ \left< U_{0,x}; E \right> - \left< U_{\tau, x}; E \right> \right] }\\
\leq &- \int_0^\tau \intO{ \left[ \left< U_{t,x}; \vr s(\vr, \vc{m}, E) \right> \partial_t \tilde \vt +
\left< U_{t,x};  s(\vr, \vc{m}, E) \vc{m} \right> \cdot \Grad \tilde \vt \right] } \dt\\
&+ \int_0^\tau \intO{ \left[ \left< U_{t,x}; \vr \tvu - \vc{m} \right> \cdot \partial_t \tilde \vu  + \left< U_{t,x}; \frac{ (\vr \tvu - \vc{m}) \otimes \vc{m} }{\vr}
\right> : \Grad \tilde \vu  \right] } \dt \\ &- (\gamma - 1) \int_0^\tau \intO{ \left[ \left< U_{t,x}; E - \frac{1}{2} \frac{|\vc{m}|^2 }{\vr} \right> \Div \tilde \vu \right] }\dt \\
&+ \int_0^\tau \intO{ \left[ \left< U_{t,x}; \vr \right> \partial_t \tvt s(\tvr, \tvt) + \left< U_{t,x}; \vc{m} \right> \cdot \Grad \tvt s(\tvr, \tvt) \right] } \dt\\
&+ \int_0^\tau \intO{ \left[ \left< U_{t,x}; \tvr - \vr \right> \frac{1}{\tvr} \partial_t p (\tvr, \tvt) - \left< U_{t,x}; \vc{m} \right> \cdot \frac{1}{\tvr} \Grad p(\tvr, \tvt) \right] } \dt\\
&+ \int_0^\tau \int_{\Omega} \Grad \tvu : {\rm d} \mu_C,
\end{split}
\end{equation}
see \cite{BreFei17}, \cite{BreFei17A}.

A similar relation holds for rDMV solutions, namely
\begin{equation} \label{rrei}
\begin{split}
\Big[ \int_\Omega &\left< U_{t,x}; \mathcal{E}_\chi \left( \vr, \vc{m}, E \Big| \tilde \vr, \tilde \vt, \tilde \vu \right) \right> \ \dx \Big]_{t = 0}^{t = \tau}
+ \intO{ \left[ \left< U_{0,x}; E \right> - \left< U_{\tau, x}; E \right> \right] }\\
\leq &- \int_0^\tau \intO{ \left[ \left< U_{t,x}; \vr \chi \left( s(\vr, \vc{m}, E) \right) \right> \partial_t \tilde \vt +
\left< U_{t,x};  \chi \left( s(\vr, \vc{m}, E) \right) \vc{m} \right> \cdot \Grad \tilde \vt \right] }\dt\\
&+ \int_0^\tau \intO{ \left[ \left< U_{t,x}; \vr \tvu - \vc{m} \right> \cdot \partial_t \tilde \vu  + \left< U_{t,x}; \frac{ (\vr \tvu - \vc{m}) \otimes \vc{m} }{\vr}
\right> : \Grad \tilde \vu  \right] } \dt \\ &- (\gamma - 1) \int_0^\tau \intO{ \left[ \left< U_{t,x}; E - \frac{1}{2} \frac{|\vc{m}|^2 }{\vr} \right> \Div \tilde \vu \right] }\dt \\
&+ \int_0^\tau \intO{ \left[ \left< U_{t,x}; \vr \right> \partial_t \tvt s(\tvr, \tvt) + \left< U_{t,x}; \vc{m} \right> \cdot \Grad \tvt s(\tvr, \tvt) \right] } \dt\\
&+ \int_0^\tau \intO{ \left[ \left< U_{t,x}; \tvr - \vr \right> \frac{1}{\tvr} \partial_t p (\tvr, \tvt) - \left< U_{t,x}; \vc{m} \right> \cdot \frac{1}{\tvr} \Grad p(\tvr, \tvt) \right] } \dt\\
&+ \int_0^\tau \int_{\Omega} \Grad \tvu : {\rm d} \mu_C,
\end{split}
\end{equation}
where
\[
\mathcal{E}_\chi \left( \vr, \vc{m}, E \Big| \tilde \vr, \tilde \vt, \tilde \vu \right)  \equiv
E - \tilde \vt \mathcal{S}_\chi(\vr, \vc{m}, E) - \vc{m} \cdot \tilde \vu + \frac{1}{2} \vr |\tilde \vu|^2 +
p(\tilde \vr, \tilde \vt) - \left( e(\tilde \vr, \tilde \vt) - \tilde \vt s( \tilde \vr, \tilde \vt) + \frac{p(\tilde \vr, \tilde \vt)}{\tilde \vr} \right) \vr.
\]

We point out that relations (\ref{rei}), (\ref{rrei}) hold for \emph{any} trio of differentiable functions $(\tvr, \tvt, \tvu)$, $\tvr, \tvt > 0$
whenever $p = p(\tvr, \tvt)$, $s = s(\tvr, \tvt)$ and $e= e(\tvr, \tvt)$ satisfy Gibbs' equation (\ref{r3}). In particular, the total entropy $\tvr s(\tvr,\tvt)$
need not be directly related to $\mathcal{S}$.

\subsubsection{DMV--strong uniqueness}

As a corollary of the relative energy inequality, we show the weak--strong uniqueness principle in the class of DMV solutions. A similar result can be obtained for the rDMV solutions, see \cite{BreFei17}. We suppose that the Euler system (\ref{i1}--\ref{i5}) possesses a smooth (continuously differentiable) solution starting from the initial data $\tvr_0$, $\tvt_0$, $\tvu_0$. In view of the specific form of the relative energy $\mathcal{E}$,
in particular its dependence on the temperature, it is convenient to express this
smooth solution in the standard variables as $(\tvr, \tvt, \tvu)$. Consequently, our goal is to show that
\[
U_{t,x} = \delta_{\tvr(t,x), \tvr \tvu(t,x), \frac{1}{2}\tvr |\tvu|^2 + \tvr e(\tvr, \tvt) (t,x) }
\ \mbox{for a.a.}\ (t,x) \in (0,T) \times \Omega
\]
for any DMV solution ${U}$ emanating from the same initial data, meaning
\[
U_{0,x} = \delta_{\tvr_0(x), \tvr_0 \tvu_0(x), \frac{1}{2}\tvr_0 |\tvu_0|^2 + \tvr_0 e(\tvr_0, \tvt_0) (x) }
\ \mbox{for a.a.}\ x \in \Omega.
\]
To this end, we substitute the solution $(\tvr, \tvt, \tvu)$ in the relative energy inequality (\ref{rei}),
with the relevant thermodynamic functions $p$, $e$, and $s$.
Our goal is to show that
\begin{equation*} %\label{EE1}
\mathcal{R}(\tau) \equiv \int_\Omega  \left< U_{\tau,x}; \mathcal{E} \left( \vr, \vc{m}, E \Big| \tilde \vr, \tilde \vt, \tilde \vu \right) \right> \ \dx
+ \intO{ \left[ \left< U_{0,x}; E \right> - \left< U_{\tau, x}; E \right> \right] } = 0
\ \mbox{for a.a.}\ \tau > 0
\end{equation*}
by applying a Gronwall type argument.
Note that $\mathcal{R}(0) = 0$.
We proceed in several steps.

\medskip

\noindent {\bf Step 1:}

By virtue of (\ref{p8}), we get
\begin{equation} \label{EE2}
\begin{split}
\mathcal{R}(\tau)
\aleq &- \int_0^\tau \intO{ \left[ \left< U_{t,x}; \vr s(\vr, \vc{m}, E) \right> \partial_t \tilde \vt +
\left< U_{t,x};  s(\vr, \vc{m}, E) \vc{m} \right> \cdot \Grad \tilde \vt \right] }\\
&+ \int_0^\tau \intO{ \left[ \left< U_{t,x}; \vr \tvu - \vc{m} \right> \cdot \partial_t \tilde \vu  + \left< U_{t,x}; \frac{ (\vr \tvu - \vc{m}) \otimes \vc{m} }{\vr}
\right> : \Grad \tilde \vu  \right] } \dt \\ &- (\gamma - 1) \int_0^\tau \intO{ \left[ \left< U_{t,x}; E - \frac{1}{2} \frac{|\vc{m}|^2 }{\vr} \right> \Div \tilde \vu \right] }\dt \\
&+ \int_0^\tau \intO{ \left[ \left< U_{t,x}; \vr \right> \partial_t \tvt s(\tvr, \tvt) + \left< U_{t,x}; \vc{m} \right> \cdot \Grad \tvt s(\tvr, \tvt) \right] } \dt\\
&+ \int_0^\tau \intO{ \left[ \left< U_{t,x}; \tvr - \vr \right> \frac{1}{\tvr} \partial_t p (\tvr, \tvt) - \left< U_{t,x}; \vc{m} \right> \cdot \frac{1}{\tvr} \Grad p(\tvr, \tvt) \right] } \dt\\
&+ \int_0^\tau \mathcal{R}(t) \ \dt,
\end{split}
\end{equation}
where $\aleq$ hides a multiplicative constant depending only on the strong solution and structural properties of the involved nonlinearities as the case may be.

\medskip

\noindent {\bf Step 2:}

Writing
\[
\left< U_{t,x}; \frac{ (\vr \tvu - \vc{m}) \otimes \vc{m} }{\vr}
\right> : \Grad \tilde \vu  = \left< U_{t,x}; \frac{ (\vr \tvu - \vc{m}) \otimes ( \vc{m} - \vr \tilde \vu) }{\vr}
\right> : \Grad \tilde \vu  + \left< U_{t,x};  \vr \tvu - \vc{m} \right> \tilde \vu \cdot \Grad \tilde \vu
\]
we deduce from (\ref{EE2}) and the fact that
\[
\partial_t \tvu + \tvu \cdot \Grad \tvu +\frac{1}{\tvr} \Grad p(\tvr, \tvt) =0,
\]
\begin{equation} \label{EE3}
\begin{split}
\mathcal{R}(\tau)
\aleq &- \int_0^\tau \intO{ \left[ \left< U_{t,x}; \vr s(\vr, \vc{m}, E) \right> \partial_t \tilde \vt +
\left< U_{t,x};  s(\vr, \vc{m}, E) \vc{m} \right> \cdot \Grad \tilde \vt \right] }\\
&- \int_0^\tau \intO{ \left< U_{t,x}; \vr \tvu \right> \cdot \frac{1}{\tvr} \Grad p(\tvr, \tvt) } \dt \\ &- (\gamma - 1) \int_0^\tau \intO{ \left[ \left< U_{t,x}; E - \frac{1}{2} \frac{|\vc{m}|^2 }{\vr} \right> \Div \tilde \vu \right] }\dt \\
&+ \int_0^\tau \intO{ \left[ \left< U_{t,x}; \vr \right> \partial_t \tvt s(\tvr, \tvt) + \left< U_{t,x}; \vc{m} \right> \cdot \Grad \tvt s(\tvr, \tvt) \right] } \dt\\
&+ \int_0^\tau \intO{ \left< U_{t,x}; \tvr - \vr \right> \frac{1}{\tvr} \partial_t p (\tvr, \tvt)  } \dt + \int_0^\tau \mathcal{R}(t) \ \dt.
\end{split}
\end{equation}

\medskip

\noindent {\bf Step 3:}

Introducing the conservative variables
\[
\tilde {m} = \tvr \tvu,\
\tilde E = \frac{1}{2} \tvr |\tvu|^2 + \frac{1}{\gamma - 1} p (\tvr, \tvt)
\]
we get
\[
\begin{split}
&\int_0^\tau \intO{ \left[ \left< U_{t,x}; \vr \right> \partial_t \tvt s(\tvr, \tvt) + \left< U_{t,x}; \vc{m} \right> \cdot \Grad \tvt s(\tvr, \tvt) \right] } \dt\\
&- \int_0^\tau \intO{ \left[ \left< U_{t,x}; \vr s(\vr, \vc{m}, E) \right> \partial_t \tilde \vt +
\left< U_{t,x};  s(\vr, \vc{m}, E) \vc{m} \right> \cdot \Grad \tilde \vt \right] }\\
&= - \int_0^\tau \intO{ \left[ \left< U_{t,x}; \tvr \nabla_{\vr, \vc{m}, E} s(\tvr, \tvm, \tvE) \left( \vr - \tvr,  \vm - \tvm,  E - \tvE \right) \right> \partial_t \tvt \right] } \dt\\
&- \int_0^\tau \intO{ \left[ \left< U_{t,x}; \tvm \nabla_{\vr, \vc{m}, E} s(\tvr, \tvm, \tvE) \left( \vr - \tvr,  \vm - \tvm,  E - \tvE \right) \right> \cdot \Grad \tvt \right] } \dt \\&+ \mbox{``quadratic'' terms}
\end{split}
\]
where the ``quadratic'' terms are controlled by $\mathcal{R}$, see \cite[Section 3.2.2]{BreFei17} for details. Accordingly, we may infer from \eqref{EE3} that
\begin{equation} \label{EE4}
\begin{split}
\mathcal{R}(\tau)
\aleq &- \int_0^\tau \intO{ \left[ \left< U_{t,x}; \tvr \nabla_{\vr, \vc{m}, E} s(\tvr, \tvm, \tvE) \left( \vr - \tvr,  \vm - \tvm,  E - \tvE \right) \right> \partial_t \tvt \right] } \dt\\
&- \int_0^\tau \intO{ \left[ \left< U_{t,x}; \tvm \nabla_{\vr, \vc{m}, E} s(\tvr, \tvm, \tvE) \left( \vr - \tvr,  \vm - \tvm,  E - \tvE \right) \right> \cdot \Grad \tvt \right] } \dt \\
&+ \int_0^\tau \intO{ \left< U_{t,x}; \tvr - \vr \right> \frac{1}{\tvr} \left(\partial_t p (\tvr, \tvt) + \tvu\cdot \Grad p(\tvr,\tvt)\right) } \dt \\
&- \int_0^\tau \intO{ \left< U_{t,x}; \tvr \right>\frac{1}{\tvr} \tvu \cdot  \Grad p(\tvr, \tvt) } \dt \\ &- (\gamma - 1) \int_0^\tau \intO{ \left[ \left< U_{t,x}; E - \frac{1}{2} \frac{|\vc{m}|^2 }{\vr} \right> \Div \tilde \vu \right] }\dt
 + \int_0^\tau \mathcal{R}(t) \ \dt.
\end{split}
\end{equation}

\medskip

\noindent{\bf Step 4:}

Next, since
\[-\int_0^\tau \intO{ \tvu \cdot \Grad p(\tvr, \tvt)}\dt = \int_0^\tau \intO{ p(\tvr, \tvt) \Div \tvu}\dt\]
and
\[
- \Div \tvu = \frac{1}{\tvr} \left( \partial_t \tvr + \tvu \cdot \Grad \tvr \right)
\]
we can rewrite (\ref{EE4}) as
\begin{equation} \label{EE6}
\begin{split}
\mathcal{R}(\tau)
\aleq &- \int_0^\tau \intO{ \left[ \left< U_{t,x}; \tvr \nabla_{\vr, \vc{m}, E} s(\tvr, \tvm, \tvE) \left( \vr - \tvr,  \vm - \tvm,  E - \tvE \right) \right> \left(  \partial_t \tvt + \tvu \cdot \Grad \tvt \right) \right] } \dt\\
&+ \int_0^\tau \intO{ \frac{1}{\tvr} \left< U_{t,x}; \tvr - \vr \right> \left( \partial_t p(\tvr, \tvt)
+ \tvu \cdot \Grad p(\tvr, \tvt) \right)} \dt \\
&+ (\gamma - 1) \int_0^\tau \intO{ \left[ \frac{1}{\tvr} \left< U_{t,x};\left( E - \frac{1}{2} \frac{|\vc{m}|^2 }{\vr} \right)
- \left( \tvE - \frac{1}{2} \frac{|\tvm|^2 }{\tvr} \right)  \right> \left( \partial_t \tvr + \tvu \cdot \Grad \tvr    \right) \right] }\dt\\
& + \int_0^\tau \mathcal{R}(t) \ \dt.
\end{split}
\end{equation}

\medskip

\noindent
{\bf Step 5:}

As shown in Section \ref{RUCV} below,
\[
\begin{split}
\tvr \partial_\vr s(\tvr, \tvm, \tvE) &= - \frac{\gamma}{\gamma - 1} \frac{ p(\tvr, \tvt) }{\tvr \tvt} + \frac{1}{2 \tvt}  |\tvu|^2, \\
\tvr \nabla_{\vc{m}} s(\tvr, \tvm, \tvE) &= - \frac{1}{\tvr \tvt} \tvm, \\
\tvr \partial_E s(\tvr, \tvm, \tvE) &= \frac{1}{\tvt}.
\end{split}
\]
Consequently, we deduce from (\ref{EE6})
\begin{equation*} %\label{EE7}
\begin{split}
\mathcal{R}(\tau)
\aleq&
\int_0^\tau \intO{ \left< U_{t,x}; \tvE - E \right> \left[ \partial_t \log \left( \frac{ \tvt }{\tvr^{\gamma - 1}} \right) +
\tvu \cdot \Grad  \log \left( \frac{ \tvt }{\tvr^{\gamma - 1}} \right) \right] } \dt
\\
&+ \int_0^\tau \intO{ \frac{\tvm}{\tvr} \cdot \left< U_{t,x}; \tvm - \vc{m} \right> \left[ \partial_t \log \left( \frac{ \tvt }{\tvr^{\gamma - 1}} \right) +
\tvu \cdot \Grad  \log \left( \frac{ \tvt }{\tvr^{\gamma - 1}} \right) \right] } \dt
\\
&+
\int_0^\tau \intO{ \frac{|\tvu|^2}{2} \cdot \left< U_{t,x}; \tvr - \vr \right> \left[ \partial_t \log \left( \frac{ \tvt }{\tvr^{\gamma - 1}} \right) +
\tvu \cdot \Grad  \log \left( \frac{ \tvt }{\tvr^{\gamma - 1}} \right) \right] } \dt
\\
& + \mbox{``quadratic'' terms} + \int_0^\tau \mathcal{R}(t) \ \dt.
\end{split}
\end{equation*}
Seeing that the entropy $s = s(\tvr, \tvt)$ is transported for the smooth solution $(\tvr, \tvt, \tvu)$
and, in accordance with (\ref{FIN}), $s$ takes the form $S(Z)$, $Z = \vt / \vr^{\gamma - 1}$, we deduce that
\[
\partial_t \log \left( \frac{ \tvt }{\tvr^{\gamma - 1}} \right) +
\tvu \cdot \Grad  \log \left( \frac{ \tvt }{\tvr^{\gamma - 1}} \right) = 0.
\]
We have shown the following result, cf. \cite{BreFei17}:

\begin{Proposition} \label{P2}

Let the thermodynamic functions $p$, $e$, and $s$ satisfy the hypotheses (\ref{r9}--\ref{r11}). Suppose that the Euler system (\ref{i1}) admits a continuously differentiable
solution $(\tvr, \tvt, \tvu)$ in $[0,T] \times \Omega$ emanating from the initial data
\[
\tvr_0 > 0,\ \tvt_0 > 0 \ \mbox{in}\ {\Omega}.
\]

Assume that  $\{ U_{t,x} \}_{(t,x) \in (0,T) \times \Omega}$ is a DMV solution of the  system (\ref{p1}), (\ref{p2}) in the sense specified in Definition \ref{D1}, such that
\[
U_{0,x} = \delta_{\tvr_0(x), \tvr_0 \tvu_0 (x), \frac{1}{2} \tvr_0(x) |\tvu_0(x)|^2 + \tvr_0 e(\tvr_0, \tvt_0)(x)} \ \mbox{for a.a.}\ x \in \Omega.
\]

Then
\[
U_{t,x} = \delta_{\tvr(t,x), \tvr \tvu (t,x), \frac{1}{2} \tvr(x) |\tvu(x)|^2 + \tvr e(\tvr, \tvt)(t,x)} \ \mbox{for a.a.}\ (t,x) \in (0,T) \times \Omega.
\]

\end{Proposition}

\section{Maximal dissipation, main result}
\label{M}

Let us start by discussing the DMV solutions satisfying the entropy inequality (\ref{p7}) with a single entropy $\mathcal{S}$. In view of the Riesz representation theorem, there exists a non--negative Borel measure $\sigma$ supported by the physical space $[0,T] \times \Omega$ such that
\begin{equation} \label{M1}
\begin{split}
\int_0^T &\intO{ \left[ \left< U_{t,x} ; \mathcal{S}(\vr, \vc{m}, E) \right> \partial_t \varphi +
\left< U_{t,x}; \mathcal{S} (\vr, \vc{m}, E) \frac{\vc{m}}{\vr} \right> \cdot \Grad \varphi  \right] } \dt + \int_0^T \int_{\Omega} \varphi \ {\rm d} \sigma  \\
&= - \intO{ \left< U_{0,x}; \mathcal{S}(\vr, \vc{m}, E) \right> \varphi(0, \cdot) }
\end{split}
\end{equation}
for any $\varphi \in \DC([0,T) \times \Omega)$, $\varphi \geq 0$.

\begin{Definition} \label{D3}
Let $U^1$, $U^2$ be two DMV solutions of the Euler system starting from the same initial data $U_0$ and satisfying (\ref{M1}) with $\sigma_1$, $\sigma_2$, respectively. We say that
\[
U^1 \succeq U^2 \ \mbox{iff} \ \sigma_1 \geq \sigma_2 \ \mbox{in}\ [0,T) \times \Omega,
\]
equivalently,
\[
\begin{split}
\int_0^T &\intO{ \left[ \left< U^1_{t,x} ; \mathcal{S}(\vr, \vc{m}, E) \right> \partial_t \varphi +
\left< U^1_{t,x}; \mathcal{S} (\vr, \vc{m}, E) \frac{\vc{m}}{\vr} \right> \cdot \Grad \varphi  \right] } \dt\\
& \leq \int_0^T \intO{ \left[ \left< U^2_{t,x} ; \mathcal{S}(\vr, \vc{m}, E) \right> \partial_t \varphi +
\left< U^2_{t,x}; \mathcal{S} (\vr, \vc{m}, E) \frac{\vc{m}}{\vr} \right> \cdot \Grad \varphi  \right] } \dt
\end{split}
\]
for any $\varphi \in \DC([0,T) \times \Omega)$, $\varphi \geq 0$.

We shall say that a DMV solution is \emph{maximal} if it is maximal with respect to the relation $\succeq$. More specifically, if
$(\Ov{U}, \Ov{\sigma})$ is maximal, and $(U, \sigma)$ is another solution of the same problem with $\sigma \geq \Ov{\sigma}$, then
$\sigma = \Ov{\sigma}$.

\end{Definition}

We are ready to formulate our main result - the existence of a maximal DMV solution.

\begin{Theorem} \label{T1}

Let the function $S$ satisfy (\ref{r11}) with $\Ov{p} > 0$, $\lim_{Z \to \Ov{p}+} S(Z) = 0$, and
\begin{equation} \label{M3}
S(Z) \leq C(1 + |\log(Z)|) \ \mbox{for all} \ Z \geq \Ov{p}.
\end{equation}
Let the initial datum $U_0$ be given such that
\begin{equation} \label{M4}
U_{0,x} \left\{ \vr > 0 \right\} = U_{0,x} \left\{ E - \frac{1}{2} \frac{ |\vc{m}|^2 }{\vr}  > \frac{\Ov{p}}{(\gamma-1)} \vr^\gamma\right\} = 1 \ \mbox{for a.a.}\ x \in \Omega,
\end{equation}
\[
\intO{ \left< U_{0,x}; E \right> } < \infty.
\]

Then the Euler system (\ref{p1}), (\ref{p2}) admits a maximal DMV solution in the sense of Definition \ref{D3}.

\end{Theorem}

\begin{Remark} \label{R4}

Hypothesis (\ref{M3}) is purely technical. It guarantees boundedness of the entropy in terms of the total energy. It is sufficient for proving the
\emph{existence} of a DMV solution, see \cite{BreFei17A}.

\end{Remark}

\begin{Remark} \label{R5}

Hypothesis (\ref{M4}) corresponds to the positivity of the initial density as well as the initial temperature.

\end{Remark}

The situation is slightly more complicated for rDMV solutions as the family of entropies is parameterized by the cut--off functions $\chi$. We introduce
a family $\{ \chi_K \}_{K = 1}^\infty$,
\[
\begin{split}
\chi_K(Z) &= K \chi \left( \frac{Z}{K} \right), \ \chi \mbox{ defined on } ( S(\Ov{p}), \infty),\ \lim_{Z \to S(\Ov{p})+} \chi(Z) = S(\Ov{p}), \\ & \chi'(Z) > 0, \ \chi''(Z) < 0,\ \chi(Z) \leq \Ov{\chi}
\ \mbox{for all}\ Z,\
\chi(Z) = Z \ \mbox{for}\ Z \leq 1.
\end{split}
\]
Now, any rDMV solution satisfies
\begin{equation} \label{M2}
\begin{split}
\int_0^T &\intO{ \left[ \left< U_{t,x} ; \mathcal{S}_{\chi_K}(\vr, \vc{m}, E) \right> \partial_t \varphi +
\left< U_{t,x}; \mathcal{S}_{\chi_K} (\vr, \vc{m}, E) \frac{\vc{m}}{\vr} \right> \cdot \Grad \varphi  \right] } \dt + \int_0^T \int_{\Omega} \varphi \ {\rm d} \sigma^K  \\
&= - \intO{ \left< U_{0,x}; \mathcal{S}_{\chi_K}(\vr, \vc{m}, E) \right> \varphi(0, \cdot) }
\end{split}
\end{equation}
for any $\varphi \in \DC([0,T) \times \Omega)$, $\varphi \geq 0$.

If the entropy complies with (\ref{M3}), and the initial data satisfy a slightly more restrictive condition than (\ref{M4}), namely
\[
U_{0,x} \left\{ \vr > \delta \right\} = U_{0,x} \left\{ E - \frac{1}{2} \frac{ |\vc{m}|^2 }{\vr}  > \frac{\Ov{p}}{(\gamma-1)} \vr^\gamma
+ \delta \right\} = 1 \ \mbox{for a.a.}\ x \in \Omega
\]
for some $\delta > 0$, it is possible to let $K \to \infty$ in (\ref{M2}) to recover (\ref{M1}). In other words, under these circumstances, a rDMV solution
is also a DMV solution and Definition \ref{D3} applies. Note that here we may allow
\[
\Ov{p} = 0 \ \mbox{as well as}\ \lim_{Z \to \Ov{p}+} S(Z) = - \infty.
\]

The remaining part of the paper is devoted to the proof of Theorem \ref{T1}. A similar statement can be formulated and proved in the context
of rDMV solutions using the above observation. We leave the details to the interested reader.

\section{Existence of maximal solutions}
\label{E}

The existence of a DMV solution under the hypotheses of Theorem \ref{T1} was established in \cite[Theorem 3.4]{BreFei17A} for $\gamma = \frac{5}{3}$.
The DMV solutions were identified as the vanishing dissipation limits of the full Navier--Stokes--Fourier system, for which the original pressure
law $p$ has been modified by adding a radiative component
\[
p_a(\vr, \vt) = p(\vr, \vt) + \frac{a}{4} \vt^4, \ a \to 0 \ \mbox{in the asymptotic limit.}
\]
A similar construction can be applied to the present case by considering
\[
p_a(\vr, \vt) = p(\vr, \vt) + a \left( \vr^{5/3} + \frac{1}{4} \vt^4 \right), \ a \to 0.
\]
Thus, under the hypotheses of Theorem \ref{T1}, the problem admits at least one DMV solution.

To establish the existence of a \emph{maximal} solution, we follow DiPerna \cite{DiP2} and use an argument based on
Zorn's lemma (Axiom of Choice). We consider the partially ordered set $\mathfrak{M}$ of all entropy production measures
(the measure satisfying (\ref{M1})) associated to DMV solutions for a given initial data $U_0$,
\[
\mathfrak{M} = \left\{ \sigma \in \mathcal{M}^+ ([0,T] \times \Omega) \ \Big| \ \sigma \ \mbox{satisfies (\ref{M1})} \right\}.
\]
Let $\mathfrak{A}$ be a chain (a totally ordered subset) in $\mathfrak{M}$. Our goal is to show that $\mathfrak{A}$ admits an upper bound
in $\mathfrak{M}$, meaning there is $\Ov{\sigma} \in \mathfrak{M}$ such that $\Ov{\sigma} \geq \sigma$ for any $\sigma \in \mathfrak{A}$.

Let $\{ g_m \}_{m=1}^\infty$ be a family of
non-negative functions in $C([0,T] \times \Omega)$, dense in the convex cone of non--negative continuous functions on the compact set
$[0,T] \times \Omega$. Denote
\[
\Ov{G}_m = \sup_{ \sigma \in \mathfrak{A} } \left< \sigma; g_m \right>.
\]
For each fixed $m$, there is a sequence $\{ \sigma_{m,n} \}_{n=1}^\infty \subset \mathfrak{A}$ such that $ \sigma_{m,1} \leq \sigma_{m,2} \leq \dots \leq \sigma_{m,n} $ and
\[\langle \sigma_{m,n};g_m\rangle
\to \overline {G}_m  \ \mbox{as}\ n \to \infty.
\]
Finally, we consider a sequence $\{ \sigma_n \}_{n=1}^\infty$,
\[
\sigma_n = \max_{m \leq n} \ \sigma_{m,n}.
\]
Thus our task reduces to finding an upper bound for the sequence of measures $\{ \sigma_n \}_{n=1}^\infty$.
Rephrased in terms of the relation $\preceq$,
we have to show
that any sequence of DMV solutions satisfying
\[
U^1 \preceq U^2 \preceq \dots \preceq  U^n
\]
admits a supremum $U \succeq U^n$ for all $n=1,2,\dots$

\subsection{Uniform bounds}

First we see that
\begin{equation} \label{E2}
\intO{ \left< U^n_{\tau, x}; E \right> } \leq \intO{ \left< U_{0, x}; E \right> } \ \mbox{uniformly in}\ n \mbox{ for a.a. }  \tau \in (0,T).
\end{equation}
Next, we deduce from the entropy balance (\ref{p7}) that
\[
\intO{ \left< U^n_{\tau,x}; \mathcal{S}(\vr, \vc{m}, E ) \right> }
\geq \intO{ \left< U_{0,x}; \mathcal{S}(\vr, \vc{m}, E ) \right> } \ \mbox{uniformly in}\ n\mbox{ for a.a. }   \tau \in (0,T),
\]
in particular
\begin{equation} \label{E3}
U^n_{t,x} \left\{ E - \frac{1}{2} \frac{ |\vc{m}|^2}{\vr} \geq \frac{\Ov{p}}{\gamma - 1} \vr^\gamma \geq 0 \right\} = 1
\ \mbox{for a.a.}\ (t,x)\mbox{ and all } n.
\end{equation}
Combining (\ref{E2}), (\ref{E3}) we may infer that
\begin{equation} \label{E4}
{\rm ess}\sup_{\tau \in (0,T)} \intO{ \left[ \left< U^n_{\tau,x}; \frac{1}{2} \frac{ |\vc{m}|^2}\vr\right> + \left< U^n_{\tau, x} ;\vr^\gamma \right> \right]
} \leq C \ \mbox{uniformly in}\ n.
\end{equation}

Finally, we may use hypothesis (\ref{M3}) and the previously established estimates to conclude that
\begin{equation} \label{E5}
{\rm ess}\sup_{\tau \in (0,T)} \intO{ \left[ \left< U^n_{\tau,x}; (\vr s)^q \right> + \left< U^n_{t,x}; (\vr s |\vc{m}|)^q \right> \right] }
\leq C \ \mbox{uniformly in}\ n \ \mbox{for some}\ q > 1.
\end{equation}

\subsection{Limit passage}

Repeating the argument of Ball \cite{BALL2}, we may assume that
\[
U^n \to U \ \mbox{weakly-(*) in}\ L^\infty_{{\rm weak}-(*)}((0,T) \times \Omega; \mathcal{M} (\mathcal{F})),
\]
where, thanks to the uniform bounds (\ref{E2}--\ref{E4}),
\[
U_{t,x} \in \mathcal{P}(\mathcal{F}) \ \mbox{for a.a.}\ (t,x).
\]
Of course, this process requires passing to a subsequence that we do not relabel here for the sake of simplicity.

Next, using again (\ref{E2}--\ref{E4}), together with (\ref{E5}), we easily observe that
\[
\int_0^T \intO{ \left[ \left< U_{t,x}; \vr \right> \partial_t \varphi + \left< U_{t,x}; \vc{m} \right> \cdot
\Grad \varphi \right] } \dt = - \intO{ \left< U_{0,x} ; \vr \right> \varphi(0, \cdot)}
\]
for any $\varphi \in \DC([0,T) \times \Omega)$, and
\[
\begin{split}
\int_0^T &\intO{ \left[ \left< U_{t,x} ; \mathcal{S}(\vr, \vc{m}, E) \right> \partial_t \varphi +
\left< U_{t,x}; \mathcal{S} (\vr, \vc{m}, E) \frac{\vc{m}}{\vr} \right> \cdot \Grad \varphi  \right] } \dt + \int_0^T \int_{\Omega} \varphi \ {\rm d} \sigma  \\
&= - \intO{ \left< U_{0,x}; \mathcal{S}(\vr, \vc{m}, E) \right> \varphi(0, \cdot) }
\end{split}
\]
for any $\varphi \in \DC([0,T) \times \Omega)$, $\varphi \geq 0$, where $\sigma \in \mathcal{M}^+([0,T) \times \Omega)$ is the weak limit of the
monotone family of measures $\sigma_1 \leq \sigma_2 \leq \dots$. In particular,
\[
\sigma \geq \sigma_n \ \mbox{for all}\ n = 1,2,\dots.
\]
Thus the proof of Theorem \ref{T1} is complete as soon as we are able to perform the limit in the momentum balance (\ref{p5}).

First, set
\[
d_n(t) = \intO{ \left< U_{0,x}, E \right> } - \intO{ \left< U^n_{t,x}; E \right> } \geq 0 \mbox{ for a.a. } t \in (0,T),
\]
and observe that, in view of (\ref{p8}),
\[
\int_0^\tau \int_\Omega {\rm d}|\mu^n_C|  \leq c(N, \gamma) \int_0^\tau d_n(t) \ {\rm d}t,
\]
where $\mu^n_C$ are the ``concentration'' measures in (\ref{p5}). Consequently, passing again to a subsequence if necessary, we may assume that
\[
\mu^n_C \to \mu_C \ \mbox{weakly-(*) in}\ \mathcal{M}([0,T) \times \Omega; R^{N \times N}),\
d_n \to d \ \mbox{weakly-(*) in}\ L^\infty(0,T),
\]
where
\[
\int_0^\tau \int_\Omega {\rm d}|\mu_C|  \leq c(N, \gamma) \int_0^\tau d(t) \ {\rm d}t.
\]

Next, thanks to the uniform bounds established in (\ref{E3}), (\ref{E4}), we may perform the limit $n \to \infty$ in (\ref{p5}) obtaining
\[
\begin{split}
\int_0^T & \intO{ \left[ \left< U_{t,x}; \vc{m} \right> \cdot \partial_t \bfphi  + \Ov{ \frac{ \vc{m} \otimes \vc{m} }{\vr}
} : \Grad \bfphi + (\gamma - 1) \Ov{ E - \frac{1}{2} \frac{|\vc{m}|^2 }{\vr} } \Div \bfphi \right] }\dt \\
&= - \intO{ \left< U_{0,x}; \vc{m} \right> \cdot \bfphi(0, \cdot) } +
\int_0^T \int_\Omega \Grad \bfphi : {\rm d} \mu_C
\end{split}
\]
for any $\bfphi \in \DC([0, T) \times \Omega; R^N)$, where
\[
\Ov{ \frac{ \vc{m} \otimes \vc{m} }{\vr}} = \mbox{weak-(*)} \lim_{n \to \infty} \left< U^n_{t,x}; \frac{\vc{m} \otimes \vc{m}}{\vr} \right>
\ \mbox{in}\ \mathcal{M}([0,T) \times \Omega; R^{N \times N}),
\]
\[
\Ov{ E - \frac{1}{2} \frac{|\vc{m}|^2 }{\vr} } = \mbox{weak-(*)} \lim_{n \to \infty} \left< U^n_{t,x};
E - \frac{1}{2} \frac{|\vc{m}|^2 }{\vr} \right> \ \mbox{in}\ \mathcal{M}([0,T) \times \Omega).
\]

Finally, we write
\[
\begin{split}
\Ov{ \frac{ \vc{m} \otimes \vc{m} }{\vr}} &= \Ov{ \frac{ \vc{m} \otimes \vc{m} }{\vr}} - \left< U_{t,x}; \frac{ \vc{m} \otimes \vc{m} }{\vr} \right>
+ \left< U_{t,x}; \frac{ \vc{m} \otimes \vc{m} }{\vr} \right>,\\
\Ov{ E - \frac{1}{2} \frac{|\vc{m}|^2 }{\vr} } &= \Ov{ E - \frac{1}{2} \frac{|\vc{m}|^2 }{\vr} }  - \left< U_{t,x}; E - \frac{1}{2} \frac{|\vc{m}|^2 }{\vr} \right> + \left< U_{t,x}; E - \frac{1}{2} \frac{|\vc{m}|^2 }{\vr} \right>.
\end{split}
\]
Note that the functions
\[
(t,x) \mapsto \left< U_{t,x}; \frac{ \vc{m} \otimes \vc{m} }{\vr} \right>,\ (t,x) \mapsto \left< U_{t,x}; E - \frac{1}{2} \frac{|\vc{m}|^2 }{\vr} \right>
\]
coincide with the so-called biting limits of the associated sequences and are integrable in $(0,T) \times \Omega$, see Ball and Murat \cite{BAMU}.
To finish the proof, we need a relation between the concentration defects
\[
\Ov{ \frac{ \vc{m} \otimes \vc{m} }{\vr}} - \left< U_{t,x}; \frac{ \vc{m} \otimes \vc{m} }{\vr} \right>,\
\Ov{ E - \frac{1}{2} \frac{|\vc{m}|^2 }{\vr} }  - \left< U_{t,x}; E - \frac{1}{2} \frac{|\vc{m}|^2 }{\vr} \right>,
\]
and the energy dissipation defect
\[
\Ov{E} - \left< U_{t,x}; E \right>.
\]
To this end, we employ the following result that can be seen as an analogue of \cite[Lemma 2.1]{FGSWW1}:
\begin{Lemma} \label{L1}
Let $\left\{ U^n_{t,x} \right\}_{n=1}^\infty$ be a sequence of parameterized probability measures on $\mathcal{F}$ such that
\[
U^n \to U \ \mbox{weakly-(*) in}\ L^\infty_{{\rm weak-} (*)}((0,T) \times \Omega; \mathcal{P}(\mathcal{F})) \ \mbox{as}\ n \to \infty,
\]
and
\[
\left\| \left< U^n_{t,x}; F(Z) \right> \right\|_{L^1((0,T) \times \Omega)} \leq C,
\]
\[
\left< U^n_{t,x}; G(Z) \right> \to \Ov{ G(Z) } , \ \left< U^n_{t,x}; F(Z) \right> \to \Ov{F(Z)} \ \mbox{weakly-(*) in} \
\mathcal{M}([0,T] \times \Omega)
\]
where $G$, $F$ are Borel functions on $\mathcal{F}$ such that
\[
|G(Z)| \leq F(Z) \ \mbox{for all}\ Z \in \mathcal{F}.
\]

Then
\[
\left| \Ov{G(Z)} - \left< U_{t,x}; G(Z) \right> \right| \leq \Ov{F(Z)} - \left< U_{t,x}; F(Z) \right>
\ \mbox{in}\ \mathcal{M}([0,T] \times \Omega).
\]

\end{Lemma}

Recalling relation (\ref{E3}), we may apply Lemma \ref{L1} to obtain the desired conclusion (\ref{p5}), (\ref{p8}).
We have shown Theorem \ref{T1}.

\section{Appendix}
\label{r}

For reader's convenience, we present some computations relating the standard variables $(\vr, \vt, \vu)$ to the conservative variables
$(\vr, \vc{m}, E)$. In particular, we clarify the relation between the thermodynamic stability hypothesis (\ref{r7}), stated in the standard variables, and concavity of the total entropy $\vr s(\vr, \vc{m}, E)$ with respect to the conservative variables, cf. also Bechtel, Rooney, and Forest \cite{BEROFO}.

\subsection{Relative energy in the standard variables}

We start by introducing the \emph{ballistic free energy}
\[
H_\Theta(\vr, \vt) = \vr \left( e(\vr, \vt) - \Theta s(\vr, \vt) \right),
\]
together with the relative energy functional
\[
\mathcal{E} \left( \vr, \vt, \vu \Big| \tilde \vr, \tilde \vt, \tilde \vu \right) =
\frac{1}{2} \vr |\vu - \tilde{\vu}|^2 + H_{\tilde \vt}(\vr, \vt) -
\frac{\partial H_{\tilde \vt}}{\partial \vr} (\tilde \vr, \tilde \vt) (\vr - \tilde \vr) - H_{\tilde \vt} (\tilde{\vr}, \tilde{\vt}),
\]
cf. \cite{FeiNov10}. Note that $\mathcal{E}$ plays the crucial role in the proof of \emph{weak--strong} and \emph{DMV--strong} uniqueness principle established in \cite{BreFei17}.

Next, we rewrite $\mathcal{E}$ as follows
\begin{equation} \label{r2}
\begin{split}
\mathcal{E} &\left( \vr, \vt, \vu \Big| \tilde \vr, \tilde \vt, \tilde \vu \right) \\ &=
\frac{1}{2} \vr |\vu|^2 + \vr e(\vr, \vt) - \vr \vu \cdot \tilde \vu + \frac{1}{2} \vr |\tilde \vu|^2 -
\tilde{\vt} \vr s(\vr, \vt) + \tilde{\vt} \tilde{\vr} s(\tilde \vr, \tilde \vt) - \tilde \vr e (\tilde \vr, \tilde \vt)\\
&- \left( e(\tilde \vr, \tilde \vt) - \tilde \vt s( \tilde \vr, \tilde \vt) \right) (\vr - \tilde \vr)
- \tilde \vr \partial_{\vr} \left( e  - \tilde \vt  s \right)(\tilde \vr, \tilde \vt) (\vr - \tilde \vr).
\end{split}
\end{equation}
Recalling Gibbs' relation (\ref{r3})
we get
\[
\tilde \vr \partial_{\vr} \left( e  - \tilde \vt  s \right)(\tilde \vr, \tilde \vt) = \frac{ p(\tilde \vr, \tilde \vt) }{\tilde \vr}.
\]
Consequently, relation (\ref{r2}) reads
\begin{equation} \label{r4}
\begin{split}
\mathcal{E} &\left( \vr, \vt, \vu \Big| \tilde \vr, \tilde \vt, \tilde \vu \right) \\ &=
\left[ \frac{1}{2} \vr |\vu|^2 + \vr e(\vr, \vt) \right] - \left[ \frac{1}{2} \tilde \vr |\tilde \vu|^2 +
\tilde \vr e(\tilde \vr, \tilde \vt) \right] + \frac{1}{2} \tilde \vr |\tilde \vu|^2  - \vr \vu \cdot \tilde \vu + \frac{1}{2} \vr |\tilde \vu|^2 \\ &-
\tilde{\vt} \left( \vr s(\vr, \vt) - \tilde{\vr} s(\tilde \vr, \tilde \vt) \right) - \left( e(\tilde \vr, \tilde \vt) - \tilde \vt s( \tilde \vr, \tilde \vt) + \frac{p(\tilde \vr, \tilde \vt)}{\tilde \vr} \right) (\vr - \tilde \vr).
\end{split}
\end{equation}
Formula (\ref{r4}) does not contain any partial derivatives of the thermodynamic functions and it is therefore easy to rewrite in the conservative variables.
This will be done in the next section.

\subsection{Relative energy in the conservative variables}
\label{RUCV}

We consider $p = p(\vr,e)$ and $s = s(\vr, e)$ as functions of the density $\vr$ and the internal energy $e$. The Gibbs relation (\ref{r3}) gives rise to
\begin{equation} \label{br1}
\frac{\partial s}{\partial e}(\vr, e) = \frac{1}{\vt}, \ \frac{\partial s}{\partial \vr}(\vr, e) = - \frac{p}{\vt \vr^2}.
\end{equation}
We  consider the conservative variables
\[
\vr, \
E = \frac{1}{2} \vr |\vu|^2 + \vr e,\ \vc{m} = \vr \vu,
\]
together with the total entropy
\[
\mathcal{S}(\vr, \vc{m}, E ) = \vr s(\vr, e) = \vr s \left( \vr, \frac{1}{\vr} \left( E - \frac{1}{2} \frac{ |\vc{m}|^2 }{\vr} \right) \right).
\]

Next, using (\ref{br1}), we compute
\[
\begin{split}
\frac{\partial \mathcal{S} (\vr, \vc{m}, E )}{\partial \vr} &= s - \frac{p}{\vt \vr} - \frac{E}{\vt \vr} + \frac{1}{\vt} \frac{|\vc{m}|^2}{\vr^2}
= \frac{1}{\vt} \left( \vt s - \frac{p}{\vr} - e + \frac{1}{2} \frac{|\vc{m}|^2}{\vr^2} \right),\\
\nabla_{\vc{m}} \mathcal{S} (\vr,  \vc{m}, E ) &= - \frac{1}{\vr \vt} \vc{m},\\
\frac{\partial \mathcal{S} (\vr,  \vc{m}, E )}{\partial E} &= \frac{1}{\vt}.
\end{split}
\]

Thus, finally, setting
\[
\tilde{E} = \frac{1}{2} \tilde \vr |\tilde \vu|^2 + \tilde \vr e(\tilde \vr, \tilde \vt), \ \tilde {\vc{m}} = \tilde \vr \tilde \vu,
\]
we deduce from (\ref{r4}) that
\begin{equation} \label{r6}
\begin{split}
\mathcal{E} &\left( \vr, \vc{m}, E \Big| \tilde \vr, \tilde \vt, \tilde \vu \right) \\ &\equiv
E - \tilde \vt \mathcal{S}(\vr,  \vc{m}, E) - \vc{m} \cdot \tilde \vu + \frac{1}{2} \vr |\tilde \vu|^2 +
p(\tilde \vr, \tilde \vt) - \left( e(\tilde \vr, \tilde \vt) - \tilde \vt s( \tilde \vr, \tilde \vt) + \frac{p(\tilde \vr, \tilde \vt)}{\tilde \vr} \right) \vr
\\
&= - \tilde \vt \left[ \mathcal{S}(\vr,  \vc{m}, E) - \mathcal{S} (\tilde \vr,  \tilde{\vm}, \tilde E) \right. \\
&\ \ \ -\left. \partial_{\vr} \mathcal{S}(\tilde \vr,  \tilde{\vm}, \tilde E)(\vr - \tilde \vr)  - \nabla_{\vc{m}} \mathcal{S}(\tilde \vr, \tilde{\vm}, \tilde E) \cdot (\vc{m} - \tilde{\vc{m}}) -
\partial_E\mathcal{S}(\tilde \vr,  \tilde{\vm}, \tilde E )(E - \tilde E)
\right].
\end{split}
\end{equation}
Equality (\ref{r6}) shows that the relative energy $\mathcal{E}$ is related to the relative entropy \` a la Dafermos \cite{Daf4} via a multiplicative factor
proportional to the absolute temperature.

\subsection{Thermodynamic stability}

Thermodynamic stability hypothesis in the standard variables reads
\begin{equation}\label{br3}
\frac{\partial p}{\partial \vr}(\vr, \vt) > 0, \ \frac{\partial e}{\partial \vt} (\vr, \vt) > 0.
\end{equation}
This is equivalent to  the statement
\begin{equation} \label{r8}
\mathcal{S}(\vr, E, \vc{m}) \ \mbox{is a concave function on its effective domain.}
\end{equation}
Indeed, it was shown in \cite{FeiNov10} that the thermodynamic stability hypothesis implies coercivity of the relative energy
\[
\mathcal{E} \left( \vr, \vt, \vu \Big| \tilde \vr, \tilde \vt, \tilde \vu \right) =
\frac{1}{2} \vr |\vu - \tilde{\vu}|^2 + H_{\tilde \vt}(\vr, \vt) -
\frac{\partial H_{\tilde \vt}}{\partial \vr} (\tilde \vr, \tilde \vt) (\vr - \tilde \vr) - H_{\tilde \vt} (\tilde{\vr}, \tilde{\vt}),
\]
specifically,
\[
\mathcal{E} \left( \vr, \vt, \vu \Big| \tilde \vr, \tilde \vt, \tilde \vu \right) \geq 0,\
\mathcal{E} \left( \vr, \vt, \vu \Big| \tilde \vr, \tilde \vt, \tilde \vu \right) = 0 \ \mbox{iff}\
(\vr, \vt, \vu) = (\tilde \vr, \tilde \vt, \tilde \vu).
\]
Consequently, (\ref{r8}) follows from (\ref{r6}).
Note that, in accordance with (\ref{br3}),
\[
\lim_{\vt \to 0+ } e(\vr, \vt) = \underline{e}(\vr) \geq 0 \ \mbox{exists for any}\ \vr > 0;
\]
in particular,
\[
E - \frac{1}{2} \frac{|\vc{m}|^2}{\vr} = \vr e > \vr \underline{e}(\vr) \geq 0\mbox{ for } \vr > 0, \ \vt > 0.
\]

\subsection{Polytropic EOS}

We suppose that
\begin{equation}\label{br5}
p = (\gamma - 1) \vr e , \ \gamma > 1.
\end{equation}
Writing $s = s(\vr,e)$ we deduce from (\ref{br1}) that
\[
\frac{\partial s(\vr,e)}{\partial \vr} = - \frac{\partial s(\vr,e)}{\partial e} \frac{p}{\vr^2} =
- \frac{\partial s(\vr,e)}{\partial e} (\gamma - 1) \frac{e}{\vr}
\]
which is a first order PDE that can be solved explicitly. We get
\[
s(\vr, e) = S \left( \frac{(\gamma - 1) e}{\vr^{\gamma - 1}} \right) =
S \left( \frac{p}{\vr^{\gamma}} \right)
\]
for a certain function $S$.

In accordance with (\ref{br1}), we have
\[
S' > 0.
\]
Moreover, the hypothesis of thermodynamic stability requires concavity of the total entropy
\[
\mathcal{S}: (\vr, \vc{m}, E) \mapsto \vr S \left( (\gamma - 1) \frac{ E - \frac{1}{2} \frac{ |\vc{m}|^2}{\vr} }{\vr^\gamma} \right)
\]
which is equivalent to the concavity of the function
\[
h: (\vr, p) \mapsto \vr S \left( \frac{ p }{\vr^\gamma} \right)\ \mbox{in the variables}\ (\vr, p).
\]
We compute
\[
\begin{split}
\frac{\partial h}{\partial \vr} = S \left( \frac{ p }{\vr^\gamma} \right) -\gamma S'\left( \frac{ p }{\vr^\gamma} \right)
\frac{p}{\vr^{\gamma}},\
\frac{\partial h}{\partial p} = \frac{1}{\vr^{\gamma - 1}} S'\left( \frac{ p }{\vr^\gamma} \right),
\end{split}
\]
and
\[
\begin{split}
\frac{\partial^2 h}{\partial \vr^2} &= - \gamma  S'\left( \frac{ p }{\vr^\gamma} \right) \frac{p}{\vr^{\gamma + 1}}
 + \gamma^2 S'\left( \frac{ p }{\vr^\gamma} \right) \frac{p}{\vr^{\gamma + 1}} + \gamma^2 S''\left( \frac{ p }{\vr^\gamma} \right) \frac{p^2}{\vr^{2\gamma + 1}},
\\
\frac{\partial^2 h}{\partial p^2} &= \frac{1}{\vr^{2 \gamma - 1}} S''\left( \frac{ p }{\vr^\gamma} \right),\\
\frac{\partial^2 h}{\partial \vr \partial p} &= \frac{1 - \gamma}{\vr^\gamma} S'\left( \frac{ p }{\vr^\gamma} \right)
- \gamma S''\left( \frac{ p }{\vr^\gamma} \right) \frac{p}{\vr^{2 \gamma}}.
\end{split}
\]

Assuming
\begin{equation} \label{rr12}
(\gamma - 1) S'(Z) + \gamma S''(Z) Z < 0 \ \mbox{for all}\ Z > 0
\end{equation}
we deduce, using also the properties (\ref{r11}) of $S$, that
\[
\frac{\partial^2 h}{\partial \vr^2} \leq 0,\ \frac{\partial^2 h}{\partial p^2} \leq 0.
\]
Moreover, we may compute the Hessian of $h$ as
\[
\begin{split}
\frac{1}{\vr^{2 \gamma}} &\left( \gamma S''(Z) \Big[ (\gamma - 1) S'(Z) Z + \gamma S''(Z) Z^2 \Big]
- \Big[ (1 - \gamma) S'(Z) - \gamma S''(Z) Z    \Big]^2 \right) \\
&= \frac{1}{\vr^{2 \gamma}} \left( - \gamma (\gamma -1 ) S''(Z) S'(Z) Z - (\gamma - 1)^2 S'(Z)^2   \right)\\
&= - \frac{(\gamma - 1)S'(Z)}{\vr^{2\gamma}} \left( \gamma S''(Z) Z + (\gamma - 1)S'(Z) \right).
\end{split}
\]
Consequently, the desired concavity of the function $h$ follows from (\ref{rr12}).

Finally, we examine the domain of definition of $S$ meaning the lower bound on the quotient $\frac{p}{\vr^\gamma}$, or, equivalently,
$\frac{e}{\vr^{\gamma - 1}}$. To this end, it is convenient to pass to the standard variables $(\vr, \vt)$. Accordingly, Gibbs' equation (\ref{r3}) gives rise
to the Maxwell relation
\[
\frac{ \partial e(\vr, \vt) }{\partial \vr} = \frac{1}{\vr^2} \left( p(\vr , \vt) - \vt \frac{\partial p (\vr, \vt)}{\partial \vt} \right);
\]
whence, by virtue of (\ref{br5}),
\begin{equation} \label{r13}
\frac{ \partial e(\vr, \vt) }{\partial \vr} = \frac{\gamma - 1}{\vr} \left(  e(\vr , \vt) - \vt \frac{\partial e (\vr, \vt)}{\partial \vt} \right).
\end{equation}
Equation (\ref{r13}) can be solved explicitly yielding $e$, or $p$, and $s$ in the form
\begin{equation} \label{FIN}
p(\vr, \vt) = (\vr \vt) \frac{ \vt^{c_v} }{\vr} P \left( \frac{\vr}{\vt^{c_v}} \right) = \frac{P(Y)}{Y^\gamma} \vr^\gamma,
\ s = s(Y),
\end{equation}
where $Y= \frac{\vr}{\vt^{c_v}}, \ c_v = \frac1{\gamma-1}$.
Moreover, as $\frac{\partial e}{\partial \vt} > 0$, we can calculate that
$\frac{P(Y)}{Y^\gamma}$ is a decreasing function of $Y$ and hence we deduce that
\[
\vr \mapsto \frac{ p(\vr, \vt) }{\vr^\gamma} \searrow \Ov{p}  \ \mbox{as}\ \vt \to 0,
\]
with some $\Ov{p} \geq 0$. The natural domain of definition of $S$ is therefore the open interval $(\Ov{p}, \infty)$.
Finally, shifting  $S$ by a constant as the case may be, we may assume that
\[
\lim_{Z \to \Ov{p}+} S(Z) \in \left\{ 0, - \infty \right\}.
\]

\def\cprime{$'$} \def\ocirc#1{\ifmmode\setbox0=\hbox{$#1$}\dimen0=\ht0
  \advance\dimen0 by1pt\rlap{\hbox to\wd0{\hss\raise\dimen0
  \hbox{\hskip.2em$\scriptscriptstyle\circ$}\hss}}#1\else {\accent"17 #1}\fi}

%\bibliographystyle{plain}
%\bibliography{citace}

\end{document}